\documentclass[
submission
]{dmtcs-episciences}

\usepackage[utf8]{inputenc}
\usepackage{subfigure}
\usepackage{epstopdf}
\usepackage{mathrsfs}
\newtheorem{theorem}{Theorem}[section]

\newtheorem{lemma}[theorem]{Lemma}
\newtheorem{corollary}[theorem]{Corollary}

\author{Koko K. Kayibi\affiliationmark{1}
\and S. Pirzada\affiliationmark{2}\thanks{Supported by SERB-DST, New Delhi under the research project number EMR/2015/001047/MS}
\and Carrie Rutherford\affiliationmark{3} }
\title[Mixing Time of Markov chain of the Knapsack Problem]{Mixing Time of Markov chain of the Knapsack Problem}
\affiliation{
  Department of Mathematics and Physics, University of Hull, UK\\
  Department of Mathematics, University of Kashmir, Srinagar, India\\
  Department of Mathematics, University of South Bank, London, UK}
\received{}
\revised{}
\accepted{}
\begin{document}
\publicationdetails{VOL}{}{ISS}{NUM}{SUBM}
\maketitle
\begin{abstract}
To find the number of assignments of zeros and ones satisfying a specific Knapsack Problem is $\#P$ hard, so only approximations are envisageable. A Markov chain allowing uniform sampling of all  possible solutions is given by Luby, Randall and Sinclair. In 2005, Morris and Sinclair, by using a flow argument, have shown that the mixing time of this Markov chain  is $\mathcal{O}(n^{9/2+\epsilon})$, for any $\epsilon > 0$. By using a canonical path argument on the distributive lattice structure of the set of solutions, we obtain an improved bound, the mixing time is given as
$$\tau_{_{x}}(\epsilon) \leq n^{3} \ln (16 \epsilon^{-1}).$$
\end{abstract}

\section{Introduction}

Let $a = (a_{_{1}}, a_{_{2}},\dots, a_{_{n}})$ be a vector in $\mathbb{R}^{^{n}}$ and  let $b$ be any real number. The \textit{Knapsack Problem}, denoted by $\mathcal{K}(n, b)$,  consists of finding  a  vector $x =(x_{_{1}}, x_{_{2}}, \dots ,x_{_{n}})$, with $x_{_{i}} \in \{0,1\}$, such that

\begin{equation} \label{Eq_Fund}
<a,x > \, = a_{_{1}}x_{_{1}}+ a_{_{2}}x_{_{2}}+\dots+a_{_{n}}x_{_{n}} \leq b.
\end{equation}

Such a vector $x$ is called  a \textit{feasible solution}. In practical terms, the problem can be phrased as follows. Assume that a traveler disposes of a knapsack which can accommodate at most $b$ kilograms. Let there be $n$ items such that item  $i$ weighs $a_{_{i}}$ kg. How can we pack up so that the total weight does not exceed the weight accomodable in the knapsack. This example indicates the relevance of the problem in various walks of life. For instance, in small to medium-sized enterprises, the  total weight  allowed may be the total investment available for a list of possible projects. For everyday travelers, the problem consists as how we can pack up our luggage so not to exceed the total weight allowed by  airlines regulations. There are many applications in business, two of these are explored in \cite{advert_appl} (marketing) and \cite{ecomm_appl} (e-commerce).

In mathematical terms, the problem consists of finding the number of solutions of the  Inequality \ref{Eq_Fund}. The vector $(0,0,\cdots,0)$  is a trivial solution. Now, it is possible to check exhaustively all the $2^{^{n}}$  possibilities and see which of these satisfy Inequality \ref{Eq_Fund}. This scheme may work for small $n$. In case $n$ is large, say $n > 100$, to check exhaustively all the possible solutions will take a long time, and by the time we are finished, the problem is obsolete (the flight is gone or the investment opportunities wasted). Hence, it is advisable to sample a solution uniformly at random.

Let $\Omega(n,b)$ be the set of all solutions of the problem  $\mathcal{K}(n,b)$. One way to sample a  solution  at random is to construct a random walk which  starts at a given element of $\Omega(n,b)$, then moves to another element according to some simple rule which changes one element to another. One such simple change consists of  changing  the entry  $x_{_{i}}$ to $x_{_{i}} + 1 $ {\emph(mod~~ 2)}.  We call such a change a \textit{flip} on the position $i$.  A  flip at a point $i$  is \textit{positive}  if $x_{_{i}}$ goes from $0$ to $1$ and a flip at a point $i$  is \textit{negative} if $x_{_{i}}$  goes from $1$ to $0$. With this terminology, we let a \textit{neutral flip} to be the fact of leaving an entry unchanged.    After running the random walk for  $t$ steps, we output  the element reached by the random walk  and consider this element, the $t^{th}$ state reached, as a random sample. Such a random walk is given in \cite{sinclair2, sinclair1} as follows.

Take the solution $x = (0,0,\dots,0)$ as the starting point. Stay on $x$ with probability $\frac{1}{2}$, and with probability {$\frac{1}{2}$ choose a position $i$ in the vector $x$ and perform a flip on the position $i$.  Call the new vector obtained by a single flip as $y$. Now, if $y$ satisfies Equation \ref{Eq_Fund}, move from $x$ to $y$, else stay at $x$.  It is routine to check that this random walk defines a Markov chain that can reach all the solutions and is aperiodic. Thus, it is an ergodic Markov chain that tends asymptotically towards a stationary state $\pi$.  So it may be used to sample uniformly a solution of $\mathcal{K}(n, b)$  problem. The key issue now is to know the time taken to reach the steady state.

Since the successive states of a Markov chain are not independent, an unbiased sample can only  be obtained  if the chain reaches stationarity, the stage when the probability of sampling a particular configuration is fixed in time. The \textit{mixing time} of a Markov chain is the number of steps necessary to reach that stationary state.  To know the mixing time of a particular Markov chain is crucial to avoid either getting a biased sample (if stationarity is not reached), or  to avoid the computational cost of running the chain more than necessary. In \cite{sinclair3}, it is shown that  the mixing time of  Markov chain defined above is polynomial on the length of the solutions of the problem. In this paper, by using a different approach which makes use of the lattice structure of the solution set $\Omega(n,b)$, we obtain a better bound.

\section{Preliminaries}

Let $\mathcal{M}$  be a Markov chain on a set of states $\Omega$. Let $P$ be the matrix of transitions from one state to another. We can visualize the Markov chain $\mathcal{M}$ as a weighted directed graph $\mathcal{G}$,  where the states are vertices of $\mathcal{G}$ and there is an edge  of weight $P(x,y)$ if the transition from state $x$ to state $y$ has probability $P(x,y)$. A Markov chain is  \textit{irreducible} if, for all pairs of states $x$ and $y$, there is an integer $t$, depending on the pair $(x,y)$, such that $P^t(x,y) > 0$. In terms of the graph $\mathcal{G}$, the Markov chain is irreducible if there is a path between every pair of vertices of $\mathcal{G}$. A Markov chain is \textit{aperiodic} if for all states $x$, there is an integer $t$ such that for all $t' \geq t $, $P^{t'}(x,x) > 0$. That is, after sufficient number of iterations, the chain has a positive probability of staying on $x$ at every subsequent step. This ensures that the return to state $x$ is not periodic. In terms of the graph $\mathcal{G}$, this can be achieved by having  a loop at every vertex.  A Markov chain is \textit{ergodic} if it is irreducible and aperiodic. If $P$ is the matrix of transitions of an ergodic Markov chain, there  is an integer $t$ such that for all the pairs of states $(x,y)$, $P^t(x,y) > 0$. (Notice that $t$ does not depend on the pair). It can also be proved that for every ergodic Markov chain with transition matrix $P$, the largest eigenvalue of $P$ is equal to 1. Using this, it can be proved that there is a unique probability vector $\pi$, such that $\pi P = \pi$. The vector $\pi$ is the \textit{stationary} distribution of $ \mathcal{M}$.  A chain  is \textit{reversible}  if $P(x,y)\pi(x) = P(y,x) \pi(y)$ for every pair of states $x$ and $y$. If a chain is reversible and the matrix $P$ is symmetric, then  $\pi(x) = \pi(y)$ for every pair $x$ and $y$. The stationary state  is then said to be \textit {uniform}. In the obvious way, we say a matrix $P$ is \textit{irreducible, ergodic, aperiodic, etc.} if the Markov chain ruled by $P$ is irreducible, ergodic, aperiodic, etc.

Let $\lambda_{_{1}} =1$ and $\lambda_{_{2}}$ be the largest and the second largest eigenvalue of the transition matrix $P$. The real number $\lambda_{_{1}}- \lambda_{_{2}}$, that is, $1- \lambda_{_{2}}$  is the \textit{spectral gap}  of the matrix $P$. It can be shown  that  the bigger the gap, the faster the mixing time of the chain. The analysis of the mixing time of a Markov chain is based on the intuition that  a random walk on the graph  $\mathcal{G}$  mixes fast ( i.e., reaches all the states quickly) if $\mathcal{G}$ has no bottleneck. That is, there is no cut between any set of vertices $S$ and its complement which blocks the flow of the Markov chain and thus prevents the Markov chain from reaching easily to some states. See \cite{jerrum1, koko1, sinclair2, sinclair1} for a better exposition on the topic, and for definitions in graph theory we refer to \cite{sp}.

Denoting the probability of $x$ at stationarity by $\pi(x)$ and the probability of moving from $x$ to $y$ by $P(x,y)$, the \textit{capacity} of the arc $e=(x,y)$, denoted by $c(e)$, is given as $c(e)=\pi(x)P(x,y)$. Let $\mathcal{P}_{_{x,y}}$ denote the set of all simple paths $p$ from $x$ to $y$ (paths that contain every vertex at most once). A \textit{flow} in $\mathcal{G}$ is a function $\phi$,  from the set of simple paths to the reals, such that
$\sum_{_{p \in \mathcal{P}_{_{x,y}}}} \phi(p)=\pi(x)\pi(y)$ for all vertices $x, y$ of $\mathcal{G}$  with $x \not =y$. An \textit{arc-flow} along an arc $e$, denoted by $\phi^{\prime}(e)$, is then defined as $\phi^{\prime}(e)=\sum_{_{p \ni e}} \phi(p)$. For a flow $\phi$, a measure of existence of an overload along an arc is given by the quantity $\rho(e)$,  where $\rho(e)=\frac{\phi^{\prime}(e)}{c(e)}$, and the \textit{cost} of the flow  $\phi$, denoted by $\rho(\phi)$, is given by $\rho(\phi)=\max_{_{e}}\rho(e)$.

If a network $\mathcal{G}$ representing a Markov chain  can support a flow of low cost, then it cannot have any bottlenecks, and hence its mixing time  should be small. This intuition is confirmed by the following theorem \cite{sinclair1}.\\

\begin{theorem}\label{Theorem 1}\cite{sinclair1} 
 Let $\mathcal{M}$ be
an ergodic reversible Markov chain  with holding probabilities
$P(x,x) \geq \frac{1}{2}$  at all states $x$. The mixing time of
$\mathcal{M}$   satisfies
$$\tau_{_{x}}(\epsilon) \leq  \rho(\phi)
|p|\bigg(ln\frac{1}{\pi(x)}+ln \frac{1}{\epsilon} \bigg),$$
where $|p|$ is the length of a longest path carrying  non-zero flow in $\phi$.
\end{theorem}

Thus one way to prove that a Markov chain mixes fast is to produce a flow along some paths where the paths and the maximal overload on edges are  polynomials on the size of the problem.

\section{Canonical path and Mixing time}

Let $\mathcal{G}$ denote the graph whose vertices are all the solutions of the $\mathcal{K}(n,b)$, and there is an edge connecting a vertex $x$ to a vertex $y$ if $x$ can be turned into $y$ by changing a single entry of the  solution  vector $x$. So $\mathcal{G}$ is just an Hasse diagram, denoted as $\mathcal{L}(n,b)$,  which is  a top-truncated lattice. For illustration, let $a_{_{1}}=5$, $a_{_{2}}= 3$, $a_{_{3}}=2$, $a_{_{4}}= 1$, and $b = 9$.  The set of solutions of the inequality $5 x_{_{1}} + 3x_{_{2}}+ 2x_{_{3}} + x_{_{4}} \leq 9$ is given in Figure \ref{Lattice(4,9)}.

\begin{figure}[h]
\centering
\includegraphics[scale=0.40]{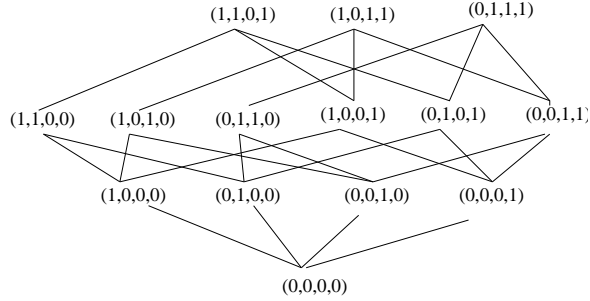}
\caption{ The solutions of the inequality $5 x_{_{1}} + 3x_{_{2}}+ 2x_{_{3}} + x_{_{4}} \leq 9$ ordered by containence.}
\label{Lattice(4,9)} \end{figure}

   Now, we  define  formally a random walk on $\mathcal{G}$ as follows. Start at any vertex $z$. If the walk is at $x$, move from $x$ to $y$ with probability $P_{_{x,y}}= \frac{1}{2n}$ if there is an edge connecting $x$ to $y$ in $\mathcal{G}$ and $<y,a> \,\, \leq b$, and stay at $x$  otherwise. Obviously, the probability of staying at $x$ is at least $\frac{1}{2}$.  Now, since $\mathcal{G}$  is just the lattice of flats of the free matroid $U_{_{n,n}}$  minus  the upper  part which contains the points  $y$ such that  $<y,a> \,\, \not \leq  b$, so $\mathcal{G}$ is connected. Therefore the matrix $P$ is  irreducible. Moreover, since  $P(x,x) \geq \frac{1}{2}$, there is a loop on every point $x$,  the matrix $P$ is aperiodic. So the random walk is an ergodic Markov chain. Finally, $P$ is a symmetric matrix. Thus, the Markov chain has a stationary uniform distribution $\pi$. We denote this Markov chain as $\mathcal{M}_{_{knap}}$.\\

The following theorem is our main result.\\

\begin{theorem}\label{fastMixing}
Let $\mathcal{L}(n,b)$ be the Hasse diagram  of the solutions of the inequality  $<a,x> \,\, \leq b$, for any real vector $a$ of length $n$. The Markov chain $\mathcal{M}_{_{knap}}$  on $\mathcal{L}(n,b)$ mixes fast and the mixing time $\tau$, is given by
$$\tau_{_{x}}(\epsilon) \leq n^{3} \ln (16 \epsilon^{-1}).$$
\end{theorem}

The proof is similar to that of papers \cite{koko1, koko2}, and requires  the following background.  A \textit{path}  from $v$  to $w$  is a sequence of binary vectors $(x^{(1)}, x^{(2)},\dots,x^{(r)})$, such that $x^{(1)} = v$ and $x^{(r)}= w$, and $x^{(i+1)}$ is obtained from $x^{(i)}$ by a single flip.
We prove this result by showing that, between any two vertices of $\mathcal{L}(n,b)$,  there is a path, the \textit{canonical path} (defined below), that is not `too long' and where  no edge is overloaded by the `probability fluid'. To put this in a more formal setting, we need the following definitions and lemmas.

Let $x_{_{i}}$ denote the $i$th entry   of $x$. Let $x$ and $y$ be two solutions of the inequality $<a,z> \,\, \leq b$. We say that $x$ and $y$ are \textit{matched} at the $i$th entry if $x_{_{i}}=y_{_{i}}$.  If $x$ and $y$ are not matched at the entry $k$, to \textit{match} the point $k$  means to do the operation $y_{_{k}} = x_{_{k}} +1 (mod\,\,  2)$. That is, if $x_{_{k}} = 0$,  change it into $1$, and vice-versa.

The \textit{canonical path} from $v$ to $w$ is the path that consists of  matching the entries in increasing order of the indices, starting from $i =1$  up to $i =n$. But, if flipping the position $i$ positively (changing its entry from $0$ to $1$) leads to a vector that is not in $\mathcal{L}(n,b)$, then we first flip negatively (change $1$ to $0$) the nearest  positions $k_{_{1}}, k_{_{2}},\dots,k_{_{r}}$, such that  $k_{_{1}} < k_{_{2}} < \dots < k_{_{r}}$,  $k_{_{j}} > i$, for $1 \leq j \leq r$,  $v_{_{k_{_{j}}}}=1$ and $w_{_{k_{_{j}}}}=0$, and $\sum_{_{j}} a_{_{k_{_{j}}}} \geq a_{_{i}}$. That is, the canonical path first  flips as many nearest positions  greater than $i$ as possible, so that flipping afterwards  the position $i$ does not lead outside  $\mathcal{L}(n,b)$. In simpler words, the Canonical path is given by Algorithm 1.

{ \bf Algorithm 1.}  To get from one feasible solution $v = (v_{_{1}}, \cdots,v_{_{n}})$ to another $w = (w_{_{1}}, \cdots,w_{_{n}})$, one scans from left to right, flipping variables from 0 to 1 and vice-versa, as required. Thus an intermediate state is a feasible solution  $(w_{_{1}}, \cdots,w_{_{i-1}},v_{_{i}},\cdots,v_{_{n}})$. However, changing $v_{_{i}} = 0$ to $w_{_{i}} = 1$ might violate the linear inequality. In this case,  one jumps $v_{_{i}}$, and  continues scanning right, changing 1s to 0s, as necessary, until there is enough slack to allow $v_{_{i}}$ to be flipped.

(Informally, one may visualize the canonical path as follows. Suppose that the items one wants to pack up in a bag are ranked from $1$ to $n$, according to their importance. Let there be  two bags, $v$ and $w$, both satisfying the requirement of not being overloaded. Assume that one would like to make $v$ carry the same items as $w$. Then  starting from the most valuable item, say item 1,  one checks whether or not the item $i$ is contained in the bag $w$. If not, and $v$ also does not  contain it, then nothing to do. In case $v$ contains it, then remove this from $v$.  Suppose that $w$ contains item $i$. If $v$ also contains it, then leave it alone. If not, check whether inserting item $i$ into $v$ would overload the bag. If there is no overload, insert item $i$ into $v$. If there is an overloading, then first  remove from $v$ the items less valuable than the item $i$ and that are not in $w$. One has to remove as many items as it is necessary to avoid the item $i$ to overload the bag.)

With this definition, suppose that Algorithm 1 is moving from $v$ to $w$, and at some instant, the highest entry flipped is $v_{_{l}}$. Then the sequence of indices  $(1, 2, \dots, n )$ can be partitioned in three parts. Part one, called the \textit{Matched Zone}, is the ordered sequence  of indices $(1,2,\dots,k )$ such that $v_{_{i}}=w_{_{i}}$, for all $i$  from $1$ to $k$. Part two, called the \textit{Heap Zone}, is the sequence of indices $(k+1, k+2, \dots, l)$, where $l$ is the  index of the last entry flipped along the canonical path, and where some entries are matched (flipped from 1 to 0  to avoid running into a binary vector that does not satisfy the inequality), and where some other entries are  not  flipped (they are jumped by Algorithm 1, either because $v_{_{j}}=0$ or $v_{_{j}}=w_{_{j}}=1$). Finally, part three, the \textit{Untouched Zone}, consists of the sequence $(l+1, l+2, \dots, n)$, the sequence of indices of entries not yet touched by Algorithm 1.

\begin{lemma} \label{Lemma 3} For all pairs of binary vectors  $v$ and  $w$ belonging to $\mathcal{L}(n,b)$, there is a canonical path  $(x^{(0)}, x^{(1)},\dots,x^{(r)})$  such that $x^{(0)} = v$ and $x^{(r)}=w$, and for all  $ 0 \leq j \leq r$,   $x^{(j)} \in \mathcal{L}(n,b)$.
\end{lemma}
\noindent{\bf Proof.} Let $v$ and  $w$ belong to $\mathcal{L}(n,b)$  but there is no such a canonical path  from $v$ to $w$. That is, Algorithm 1 gets stalled at some point. This is possible only if there is an entry $v_{_{i}}$ whose flipping positively will lead to a vector $v^{\prime}$  such that $v^{\prime} \not \in \mathcal{L}(n,b)$, but for all $j$ such that $j >i$ , $v_{_{j}}=1, w_{_{j}} =0$, we have $a_{_{i}} > \sum_{_{j}} a_{_{j}} $. That is, even after Algorithm 1  had flipped negatively all such entries $v_{_{j}}$, flipping $v_{_{i}}$  would still lead to a binary vector    $v'$   where  $  b < \sum_{_{i}} v_{_{i}}^{\prime}.a_{_{i}} $  (That is, $v^{\prime}$ is outside $\mathcal{L}(n,b)$).

Now, suppose that the canonical  path is at its $t$th  iteration when it has to flip  the entry $v_{_{i}}$. Thus,  there are $n-t$  positions left to match before reaching $w$. Let there be $s$ positions $j$  such that  $j >i$ , $v_{_{j}}=1, w_{_{j}} =0$  (points that may be flipped  to avoid exiting $\mathcal{L}(n,b)$). ( Obviously,  $n-t \geq s$, since  there may still be  positions $k$ such that   $v_{_{k}}=1, w_{_{k}} =1$ or   $v_{_{k}}=0, w_{_{k}} =0$, or  $v_{_{k}}=0, w_{_{k}} =1$).  Since flipping all the $s$ positions $j >i$ such that $v_{_{j}}=1, w_{_{j}} =0$ still leads to a binary vector outside $\mathcal{L}(n,b)$, we have $a_{_{i}} > \sum_{_{j}} a_{_{j}}$. Thus,  $w \not \in  \mathcal{L}(n,b)$. This is clearly a contradiction. \qed \\

One of the most important properties of Algorithm 1  is that once matched, an entry can not be un-matched anymore. That is, the matched zone only increases along the canonical path.  This entails the following lemma. \\

\begin{lemma}\label{lengthOfPath}
For all binary vectors  $v$ and $w$ that belong to $\mathcal{L}(n,b)$, the length of the canonical path from $v$ to $w$  is at most $n$.
\end{lemma}
\noindent{\bf Proof.}
 Since every entry is matched at most once, the length of the path is at most the number of entries in $v$.\qed \\

In order for the canonical paths argument to work, we need the number of paths through any edge of the upper-truncated lattice $\mathcal{L}(n,b)$ to be bounded above by a small polynomial of $n$ times $N$, the number of valid knapsack solutions. And indeed, in what follows, for any arc  $e=(z,y)$ of $\mathcal{L}(n,b)$, the number of different canonical paths  passing through $e$ is at most $2N$.  But, before proving this crucial result in Corollary \ref{Corollary 8}, we need the following  definitions and lemmas.  Let $v$ and $z$ be elements of $\mathcal{L}(n,b)$.  We say an entry $v_{_{i}}$ is \textit{determined} by another entry $z_{_{i}}$ if it is possible  to know the value of   $v_{_{i}}$ from the value of $z_{_{i}}$.  That is, $v_{_{i}}$ and $z_{_{i}}$  are not independent, since either $v_{_{i}}= z_{_{i}}$, or $v_{_{i}}= z_{_{i}} +1 \textrm{(mod 2)}$. We say that the index $i$ is \textit{free} in a vector $v$ if  $v_{_{i}}=0$ or $v_{_{i}}= 1$ independently of any entry of $z$. We recall that the canonical path $(v,\dots ,w)$ passes through $e=(z,y)$  means that to change $v$ into $w$, Algorithm 1 must make a sequence of flips that change  $v$ into $z$, then $z$ into $y$, then $y$ into $w$.  See Figure \ref{atMostNPath} for an illustration. \\

\begin{figure}[h]
\center
\includegraphics[scale=0.2]{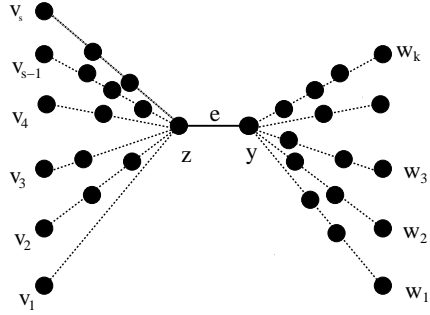}
\caption{Different canonical paths that pass through the arc $e$.}
\label{atMostNPath}
\end{figure}

\begin{lemma} \label{Lemma 5} Let $\mathcal{G}$  be the upper-truncated lattice  whose vertices are all the elements of $\mathcal{L}(n,b)$. Let $z$ and $y$ be two such vertices,  such  that $z_{_{i}} = y_{_{i}}$ for all the indices except $r$, and let  $e=(z,y)$ be the edge that represents the flip matching the $r$th entry of $z$ (That is, the flip transforming $z$ into $y$). For any two vertices $v$ and $w$, the canonical path $(v,\dots ,w)$ passes through $e$ only if the set of indices that are  free in $v$ is disjoint from the set of indices that are  free in $w$.
\end{lemma}
\noindent{\bf Proof.} Suppose that the canonical path from $v$ to $w$ passes through the arc $e$ at the $r$th entry of $z$.  Consider the Matched Zone, the Untouched Zone and the Heap Zone  at the instant when the flip is performed at $r$.  We aim at  showing the following.

\begin{itemize}
\item[(1)] The Matched Zone (from index  $1$ to $k$) is free in $v$  and determined by $y$ in $w$.
\item[(2)] The Untouched Zone (from index  $l+1$ to $n$)  is free in $w$ and determined by $z$ in $v$.
\item[(3)] In the Heap Zone, entries that  are flipped carry $0$ in $z$, and  must carry $1$ in $v$.
\item[(4)] In the Heap Zone, entries that are not flipped  and carry $0$ in $z$  must carry  $0$ in $v$.
In the Heap Zone, entries that are not flipped  and carry $1$ in $z$  must carry  $1$ in $v$.
\end{itemize}

Indeed, the canonical path from $v$ to $w$ passes through $e=(z,y)$ only if the sequence of  flips from the first one to the $l$th flip change $v$ to $z$. Thus, (1), all the entries of $v$ in the Untouched Zone  must be equal to  the corresponding entries in $z$, since there are untouched, but they are are free in $w$. Now, (2), after  changing $z$ to $y$, the remaining flips must change $y$ to $w$. Since the canonical path does not un-match positions that  are already matched, entries of  $w$ in the Matched Zone  must be equal to the corresponding entries in $y$, while  they are free  in $v$.
 In the  Heap zone, entries are either flipped  or left untouched (jumped). If, (3), they were flipped, thus they are $0$ in $z$ and they must be $1$ in $v$. Those are the entries that are flipped from $1$ to $0$ to get some slack. If, (4), they are jumped, they may be $0$ or $1$ in $z$. If they are $0$ in $z$, they must be $0$ in $v$, and if they  are $1$ in $z$, they are also $1$ in $v$, since  the Canonical Path from $v$ to $w$ reaches $z$ at that instant  means that the untouched entries are equal. Hence in the Heap zone, entries of $v$ are determined by  entries of $z$.\qed \\

\begin{lemma}\label{lemma 6} Let $\Omega(n,b)$ be the set of  solutions of the inequality $a_{_{1}}x_{_{1}}+a_{_{2}}x_{_{2}}+ \cdots + a_{_{n}}x_{_{n}}\leq b$ and let $|\Omega|=N$. Then
\begin{equation}
N= N^{\prime} + N^{\prime\prime},
\end{equation}
where $N^{\prime}$ is the number of solutions  of the inequality $a_{_{2}}x_{_{2}}+ \dots + a_{_{n}}x_{_{n}}\leq b$ and $N^{\prime\prime}$ is the number of solutions  of the inequality $a_{_{2}}x_{_{2}}+ \dots + a_{_{n}}x_{_{n}}\leq b-a_{_{1}}$.
\end{lemma}
\noindent{ \bf Proof.}
(First, note that there is nothing special with choosing  to remove $x_{_{1}}$. Any other variable would do as well.) The set of  solutions of the inequality $a_{_{1}}x_{_{1}}+a_{_{2}}x_{_{2}}+ \dots + a_{_{n}}x_{_{n}}\leq b$ can be partitioned into the sets $A$ and $B$, where $A$ is the set of solutions with $x_{_{1}}=0$ and $B$ is the set of solutions with $x_{_{1}}=1$. The set $A$ is in one-one correspondence with the set of solutions of the the inequality $a_{_{2}}x_{_{2}}+ \dots + a_{_{n}}x_{_{n}}\leq b$, and the set $B$ is in one-one correspondence with the set of solutions of inequality $a_{_{2}}x_{_{2}}+ \dots + a_{_{n}}x_{_{n}}\leq b-a_{_{1}}$.
\qed \\

\begin{corollary}\label{Corollary 7} Let $z$ be a given feasible solution, let $\alpha+\beta =n$, let  $|\Omega|=N$, and  let $V$  be the set of all the solutions where the first $\beta$  entries  are determined by  $z$. Then $|V|$ is less or equal to $(2N)^{^{\frac{\alpha}{n}}}$.
\end{corollary}
\noindent{\bf Proof.} (Obviously,  $|V|$ is just the number of solutions $a_{_{\beta+1}}x_{_{\beta+1}}+ \dots + a_{_{n}}x_{_{n}}\leq b-a_{_{1}}-a_{_{2}}-\cdots-a_{_{\beta}}$).

The proof uses induction on $\alpha$, the number of indices that are free. If $\alpha = 0$ (i.e., all the variables $x_{_{i}}$ are determined), the number of solutions is less or equal to $N^{^{0}}=1$. (We may also check the formula by taking $\alpha =n$. If all the variables are free,  there are  $N^{^{\frac{n}{n}}}=N$ possible  feasible solutions).
 For induction, assume that the result holds for $\alpha -1$.  Let $M$ denote the number of feasible solutions where $\alpha$  variables are free, and suppose that  $x_{_{n}}$ is free. Then, by Lemma \ref{lemma 6}, $M =  M^{\prime} + M^{\prime\prime}$, where $M^{\prime}$ is the number of solutions  of the inequality $a_{_{1}}x_{_{1}}+ \dots + a_{_{n}}x_{_{n-1}}\leq b$ and $M^{\prime\prime}$ is the number of solutions  of the inequality $a_{_{1}}x_{_{1}}+ \dots + a_{_{n}}x_{_{n-1}}\leq b-a_{_{n}}$. We have
\begin{eqnarray*}
M &=&  M' + M''\\
  &\leq &  N'^{^{\frac{\alpha -1}{n}}} + N''^{^{\frac{\alpha -1}{n}}} \,\,\,\,\,\textrm{by inductive hypothesis}\\
  &\leq & (2N')^{^{\frac{\alpha -1}{n}}}  \,\,\,\,\,\textrm{ since} \,\,\,\,N' \geq N''\\
  &\leq& (2N)^{^{\frac{\alpha}{n}}}.
\end{eqnarray*}

\qed \\

\begin{corollary}\label{Corollary 8}
If $e=(z,y)$ is an arc of $\mathcal{L}(n,b)$, the number of different canonical paths  passing through $e$ is at most $2N$,
  where $N$ is the number of vertices of $\mathcal{L}(n,b)$.
\end{corollary}

\noindent{\bf Proof.} Let the canonical path from $v$ to $w$ passes through the arc $e$. By Lemma \ref{Lemma 5}, the set of indices that are free in $v$ is disjoint from the set of indices that are  free in $w$. Let  $\alpha$ and $\beta$, with $\alpha + \beta \leq n$, be the number of indices that are free in $v$ and $w$, respectively. By Corollary \ref{Corollary 7}, for a fixed feasible solution $w^{\prime}$, there are at most $N^{^{\frac{\alpha}{n}}}$ vectors $v$ such that the canonical path $(v,\dots,w^{\prime})$ passes through $e=(z,y)$, and for a fixed vertex $v^{\prime}$, there are at most $N^{^{\frac{\beta}{n}}}$ vectors $w$  such that the canonical path $(v^{\prime},\dots,w)$ passes through $e=(z,y)$.   Therefore, for a fixed vector $z$ and fixed index $k$ such that $e$ is the  edge consisting of flipping the $k$th entry of $z$, the total number of canonical paths passing through  $e$  is at most $(2N)^{^{\frac{\alpha}{n}}}(2N)^{^{\frac{\beta}{n}}} \leq  2N$. \qed \\

\noindent {\bf Proof of Theorem \ref{fastMixing}.} In order to prove Theorem \ref{fastMixing} by using Theorem \ref{Theorem 1}, we  show that there is a flow $\phi$ such that $\rho(\phi)$ is a polynomial in $n$, the size of a vector solution of $\mathcal{K}(n,b)$. Indeed, if $x$ and $y$ are two vertices of $\mathcal{G}$, let $\hat{p}_{_{xy}}$ denote the canonical path from $x$ to $y$ and let $\phi$ be a flow consisting of injecting $\pi(x)\pi(y)$ units of flow along $\hat{p}_{_{xy}}$. Then, for all arcs $e$, we have $$\phi^{\prime}(e)=\sum \pi(x)\pi(y),$$ where the sum is over all the pairs $\{x,y\}$ such that $e \in \hat{p}_{_{xy}}$. Since by Corollary \ref{Corollary 8}, there are at most $2N$ canonical paths through $e$, and $\pi(x)=\pi(y)=\frac{1}{N}$, as the distribution $\pi$ is uniform, we have $$\phi^{\prime}(e) \leq  2N \pi(x)\pi(y) \leq \frac{2}{N}.$$ Moreover, since $P_{_{x,y}} = \frac{1}{2n}$ if there is an edge from $x$ to $y$, we have
\begin{equation}\label{eq2}
c(e)=\pi(x)P_{_{x,y}}\geq \frac{1}{2nN}.
\end{equation}

Thus,

\begin{equation} \label{eq3}
\rho(\phi) \leq \frac{max_{_{e}} \phi^{'}(e)}{min_{_{e}}c(e)} \leq 4 n.
\end{equation}

Now, using Theorem \ref{Theorem 1}, Lemma \ref{lengthOfPath} and  Equation \ref{eq3}, we have

$$\tau_{_{x}}(\epsilon) \leq (n)(4n)(\ln\frac{1}{\pi(x)}+\ln \frac{1}{\epsilon}).$$

Finally, using the fact that $ \frac{1}{\pi(x)}= N \leq 2^{^{n}}$, we have

$$\tau_{_{x}}(\epsilon) \leq 4n^2(n\ln 2 +\ln \frac{1}{\epsilon}) \leq n^3\ln(16\epsilon^{-1}).$$
$\Box$

{\bf Acknowledgements.} The research of second author is supported by SERB-DST, New Delhi under the research project number EMR/2015/001047/MS.

\end{document}